\documentclass[twocolumn,twoside]{IEEEtran}
\pdfoutput=1
\usepackage{amsmath,amssymb,amsthm,epsfig,dsfont,color,subfigure,empheq,graphicx,subfig}
\usepackage{enumerate,url,algpseudocode,algorithm,wasysym,epstopdf}
\usepackage[table]{xcolor}
\usepackage{arydshln}
\usepackage{tikz}
\usetikzlibrary{matrix,decorations.pathreplacing}

\DeclareMathOperator{\rank}{rank}

\DeclareMathOperator{\diag}{dg}

\DeclareMathOperator{\trace}{Tr}

\DeclareMathOperator*{\find}{find}

\newtheorem{lemma}{Lemma}
\newtheorem{corollary}{Corollary}
\newtheorem{theorem}{Theorem}

\newtheorem{criterion}{Criterion}
\theoremstyle{remark}\newtheorem{remark}{Remark}
\allowdisplaybreaks

\begin{document}
\title{Enhancing Observability in Distribution Grids\\
	using Smart Meter Data}

\author{
	Siddharth Bhela,~\IEEEmembership{Student Member,~IEEE,}
	Vassilis Kekatos,~\IEEEmembership{Member,~IEEE,} and
	Sriharsha Veeramachaneni

	\vspace*{-1.5em} 

\thanks{The work of S. Bhela was supported by a grant from Windlogics Inc.}


}


\maketitle
\begin{abstract}
Due to limited metering infrastructure, distribution grids are currently challenged by observability issues. On the other hand, smart meter data, including local voltage magnitudes and power injections, are communicated to the utility operator from grid buses with renewable generation and demand-response programs. This work employs grid data from metered buses towards inferring the underlying grid state. To this end, a coupled formulation of the power flow problem (CPF) is put forth. Exploiting the high variability of injections at metered buses, the controllability of solar inverters, and the relative time-invariance of conventional loads, the idea is to solve the non-linear power flow equations jointly over consecutive time instants. An intuitive and easily verifiable rule pertaining to the locations of metered and non-metered buses on the physical grid is shown to be a necessary and sufficient criterion for local observability in radial networks. To account for noisy smart meter readings, a coupled power system state estimation (CPSSE) problem is further developed. Both CPF and CPSSE tasks are tackled via augmented semi-definite program relaxations. The observability criterion along with the CPF and CPSSE solvers are numerically corroborated using synthetic and actual solar generation and load data on the IEEE 34-bus benchmark feeder.
\end{abstract}


\section{Introduction}\label{sec:intro}
Power flow (PF) and power system state estimation (PSSE) are central to planning, monitoring and control of electricity networks. Although these problems have been extensively studied in transmission networks~\cite{Mo00}, there is renewed interest in developing novel techniques for determining the system state in distribution grids~\cite{KlauberZhu15}. Due to limited instrumentation, low investment interest in the past, and the sheer scale of residential electricity networks, low-voltage grids have limited observability~\cite{Baran01}. Traditionally, utility operators monitor distribution grids by collecting load, voltage and current magnitude measurements only at a few buses. This mode of operation has been functional due to the stationarity of conventional loads, the availability of historical data, and the underutilization of distribution grids. Nevertheless, with the advent of distributed renewable generation, electric vehicles, and demand response, there is a need for enhancing grid observability both in space and time to accomplish grid dispatch objectives, such as voltage regulation or power loss minimization. At the same time, smart meter data are readily available at high temporal resolution; while the smart inverters found in solar panels, energy storage units, and electric vehicles can be controlled within microseconds. Both technologies could be leveraged for inferring the system state.

Determining whether a transmission network is observable given a set of specifications is performed using topological or numerical methods~\cite{Mo00}. Topological methods rely on graph-theoretic principles, whereas the numerical ones study the rank of the Jacobian matrix related to the decoupled power flow model. An analysis of PF problems under different specification buses is conducted in~\cite{Guo13}, again under the decoupled model. Given the higher resistance-to-reactance ratios, the decoupled grid model does not apply to distribution systems. 

Techniques for placing limited meters to improve distribution system state estimation (DSSE) have been reported in \cite{Baran96}, \cite{Shafiu05}, \cite{RSingh09}. A heuristic rule aiming at reducing the variance of voltage magnitude estimates at non-metered buses is suggested in \cite{Shafiu05}, and is extended to voltage angles too in \cite{RSingh09}. Pseudo-measurements have also been used in restoring observability in distribution grids~\cite{KAClem11}. Historical load profiles and real-time measurements are input to a neural network-based state estimator in~\cite{Manitsas12}, while a robust DSSE using pseudo-measurements is developed in~\cite{JYN13}. 

Beyond pseudo-measurements, smart meter and synchrophasor data have been lately utilized for various distribution inference tasks. Reference \cite{AFD15} proposes a method for adjusting the mean and variance on the noise of each smart meter based on its communication delay. Alternative weighted least-squares (WLS) and least absolute value (WLAV) solvers relying on semidefinite program (SDP) relaxations have been devised in \cite{KlauberZhu15}. Graph-theoretic tests based on the properties of the inverse covariance of voltage magnitude data have been developed in \cite{BoSch13} and \cite{Deka1} for distribution topology detection. 

The scheme of \cite{Arghandeh15} uses micro-synchrophasors to select the topology attaining the smallest state estimation error. Synchronized measurements are also employed in \cite{LiaoWengRajagopal16} for topology detection via $\ell_1$-norm penalized linear regression. The idea of collecting synchronized data to determine network topology and parameters is explored in \cite{Ardakanian17}. Nevertheless, synchrophasors are not widely employed in current low-voltage distribution grids due to deployment costs. The trade-off between state estimation accuracy and investment cost using synchrophasor and smart meter data is explored in~\cite{JLiu12}.

This work builds on smart meter data to infer the state of distribution systems. Our contribution is three-fold. Firstly, exploiting the stationarity of conventional non-metered loads and the controllability of metered buses, we put forth the idea of jointly recovering successive system states through the coupled power flow (CPF) formulation of Section~\ref{sec:problem}. Secondly, the observability of the related CPF equations in radial grids is characterized in Section~\ref{sec:analysis} via an intuitive graph-theoretic criterion connecting metered and non-metered buses on the underlying grid. The criterion was presented without proofs in the conference precursor of this work~\cite{BKVZ17}. Thirdly, the proposed CPF task is tackled in both its noiseless and noisy alternatives through SDP-based solvers presented in Section~\ref{sec:solvers}, and numerically validated using actual solar generation and residential load data from the Pecan Street project \cite{pecandata} on the IEEE 34-bus benchmark feeder in Section~\ref{sec:sim}. Conclusions and current research efforts are outlined in Section~\ref{sec:conclusions}.

\emph{Notation}: Column vectors (matrices) are denoted using lower- (upper-) case boldface letters. Sets are denoted using calligraphic symbols, while $|\mathcal{X}|$ denotes the cardinality of set $\mathcal{X}$. The notation $X_{nm}$ denotes the $(n,m)$-th entry of $\mathbf{X}$. Symbols $(\cdot)^{\top}$ and $(\cdot)^H$ stand for (complex) transposition. The operator $\diag(\mathbf{x})$ defines a diagonal matrix having $\mathbf{x}$ on its main diagonal, while $\trace(\cdot)$ and $\rank(\cdot)$ denote the matrix trace and rank, respectively.

\section{Problem Formulation}\label{sec:problem}

\subsection{Modeling Preliminaries}\label{subsec:model}
A distribution grid having $N+1$ buses can be represented by a graph $\mathcal{G}_d=(\mathcal{N}^+,\mathcal{L})$, whose nodes $\mathcal{N}^+:=\{0,\ldots,N\}$ correspond to buses and edges $\mathcal{L}$ to distribution lines. The substation bus is indexed by $n=0$; all other buses constitute set $\mathcal{N}:=\{1,\ldots,N\}$. The complex voltage at bus $n$ is expressed in polar or rectangular coordinates as $V_n=|V_n|e^{j\theta_n}=v_{r,n}+jv_{i,n}$. The $2(N+1)$--length vector collecting the rectangular components of voltages across all buses constitutes the system state $\mathbf{v}$. The latter can be partitioned as $\mathbf{v}:=[\mathbf{v}_r^\top~\mathbf{v}_i^\top]^\top$, where $\mathbf{v}_r:=[v_{r,0}~\ldots~v_{r,N}]^\top$ and $\mathbf{v}_i:=[v_{i,0}~\ldots~v_{i,N}]^\top$.

The squared voltage magnitude at bus $n$ can be expressed as a quadratic function of $\mathbf{v}$:
\begin{equation}\label{eq:PF0}
	|V_n(\mathbf{v})|^2=v_{r,n}^2 + v_{i,n}^2.
\end{equation}
To express power injections as functions of $\mathbf{v}$, let us introduce the bus admittance matrix $\mathbf{Y}:= \mathbf{G} + j\mathbf{B}$ satisfying the property $\mathbf{Y}\mathbf{1}=\mathbf{G}\mathbf{1}=\mathbf{B}\mathbf{1}=\mathbf{0}$ if shunt elements are ignored~\cite{ExpConCanBook}.  Then, the active and reactive power injections $p_n$ and $q_n$ into bus $n$ can be expressed as quadratic functions of $\mathbf{v}$ as well:
\begin{subequations} \label{eq:PF}
	\begin{align}
	p_n(\mathbf{v})= &v_{r,n} \sum_{m=1}^N \left(v_{r,m}G_{nm}-v_{i,m}B_{nm}\right) \nonumber\\
	&+v_{i,n} \sum_{m=1}^N\left(v_{i,m}G_{nm}+v_{r,m}B_{nm}\right) \label{eq:PFa}\\
	q_n(\mathbf{v})= &v_{i,n} \sum_{m=1}^N \left(v_{r,m}G_{nm}-v_{i,m}B_{nm}\right) \nonumber\\
	&-v_{r,n} \sum_{m=1}^N \left(v_{i,m}G_{nm}+v_{r,m}B_{nm}\right).  \label{eq:PFb}
	\end{align}
\end{subequations}

In the conventional PF problem, the system operator fixes two out of the three quantities $(p_n,q_n,|V_n|^2)$ for each bus $n\in\mathcal{N}$ to specified values, and solves \eqref{eq:PF0}--\eqref{eq:PF} to recover the underlying state $\mathbf{v}$. Depending on which pair of quantities is fixed, bus $n$ can be classified as:
\begin{itemize}
	\item a PQ bus if $(p_n,q_n)$ are specified;
	\item a PV bus if $(p_n,|V_n|^2)$ are specified; or
	\item the reference bus for which $|V_n|^2$ is specified and $\theta_n=0$.
\end{itemize}
In distribution grids, the substation bus is selected as the reference bus and the remaining buses are usually modeled as PQ buses~\cite{BW3}. Different from transmission grids though, there are many distribution grid buses with all three quantities $(p_n,q_n,|V_n|^2)$ unknown.

\subsection{Coupling Power Flow Problems}\label{subsec:cpf}
Since power injections in smart grids can vary unpredictably over short periods of time, power flow approaches relying on pseudomeasurements may be inadequate. On the other hand, smart meter and power inverter readings could be used to find the grid state. Different from the conventional PF setup, we classify buses into four sets: 
\begin{itemize}
\item the set $\mathcal{M}$ of metered buses for which all three quantities $(p_n, q_n,|V_n|^2)$ are specified; 
\item the set $\mathcal{O}$ of non-metered buses with no information; 
\item the singleton $\mathcal{S}=\{0\}$ of the substation; and 
\item the remaining buses in $\mathcal{N}^+$ comprising the set $\mathcal{C}$ for which $(p_n,q_n)$ are specified. 
\end{itemize}
When the distribution grid is at state $\mathbf{v}$, specifications at the aforementioned buses are provided as follows
\allowdisplaybreaks
\begin{subequations}\label{eq:cpf1}
	\begin{align}
		v_{i,0}(\mathbf{v})^2 &=0\label{eq:cpf1:0}\\
		|V_n(\mathbf{v})|^2&=|\hat{V}_n|^2 & \forall n \in \mathcal{S} \cup \mathcal{M}\label{eq:cpf1:a}\\
		q_n({\mathbf{v}})&= \hat{q}_n & \forall n \in \mathcal{C} \cup \mathcal{M}\label{eq:cpf1:c}\\
		p_n(\mathbf{v})&=\hat{p}_n & \forall n \in \mathcal{C} \cup \mathcal{M}\label{eq:cpf1:b}
	\end{align}
\end{subequations}
where hatted symbols denote specified values. Apparently, if $|\mathcal{M}|=M$ and $|\mathcal{O}|=O$, then $|\mathcal{C}|=N-M-O$. A necessary condition for solving \eqref{eq:cpf1} requires $3M+2C+2\geq 2N+2$ or $M\geq 2O$, that is the number of metered buses must be at least twice the number of non-metered ones. 

To relax the requirement $M\geq 2O$ on metered buses, a coupled version of the problem in \eqref{eq:cpf1} is considered. The idea is to couple the power flows corresponding to state $\mathbf{v}$ with the power flows corresponding to a different state $\mathbf{v}'$ observed at a subsequent time instant. By collecting the grid outputs for state $\mathbf{v}'$, the ensuing specifications are obtained:
\begin{subequations} \label{eq:cpf2}
	\begin{align} 
		v_{i,0}(\mathbf{v}')^2 &=0\label{eq:cpf2:0}\\
		|V_n(\mathbf{v}')|^{2}&=|\hat{V}_n'|^{2} & \forall n \in \mathcal{S} \cup \mathcal{M}\label{eq:cpf2:v}\\
		q_n(\mathbf{v}')&=\hat{q}_n' & \forall n \in \mathcal{C} \cup \mathcal{M}\label{eq:cpf2:q}\\
		p_n({\mathbf{v}'})&=\hat{p}_n' & \forall n \in \mathcal{C} \cup \mathcal{M}.\label{eq:cpf2:p}
	\end{align}
\end{subequations}
The equations in \eqref{eq:cpf2} exhibit a structure identical to the one in \eqref{eq:cpf1}. If the system of equations in \eqref{eq:cpf1} is under-determined, so is the system in \eqref{eq:cpf2}. Apparently, solving \eqref{eq:cpf1}--\eqref{eq:cpf2} jointly provides no added benefit towards recovering $(\mathbf{v},\mathbf{v}')$. 

\begin{figure}[t]
	\centering
	\includegraphics[scale=0.36]{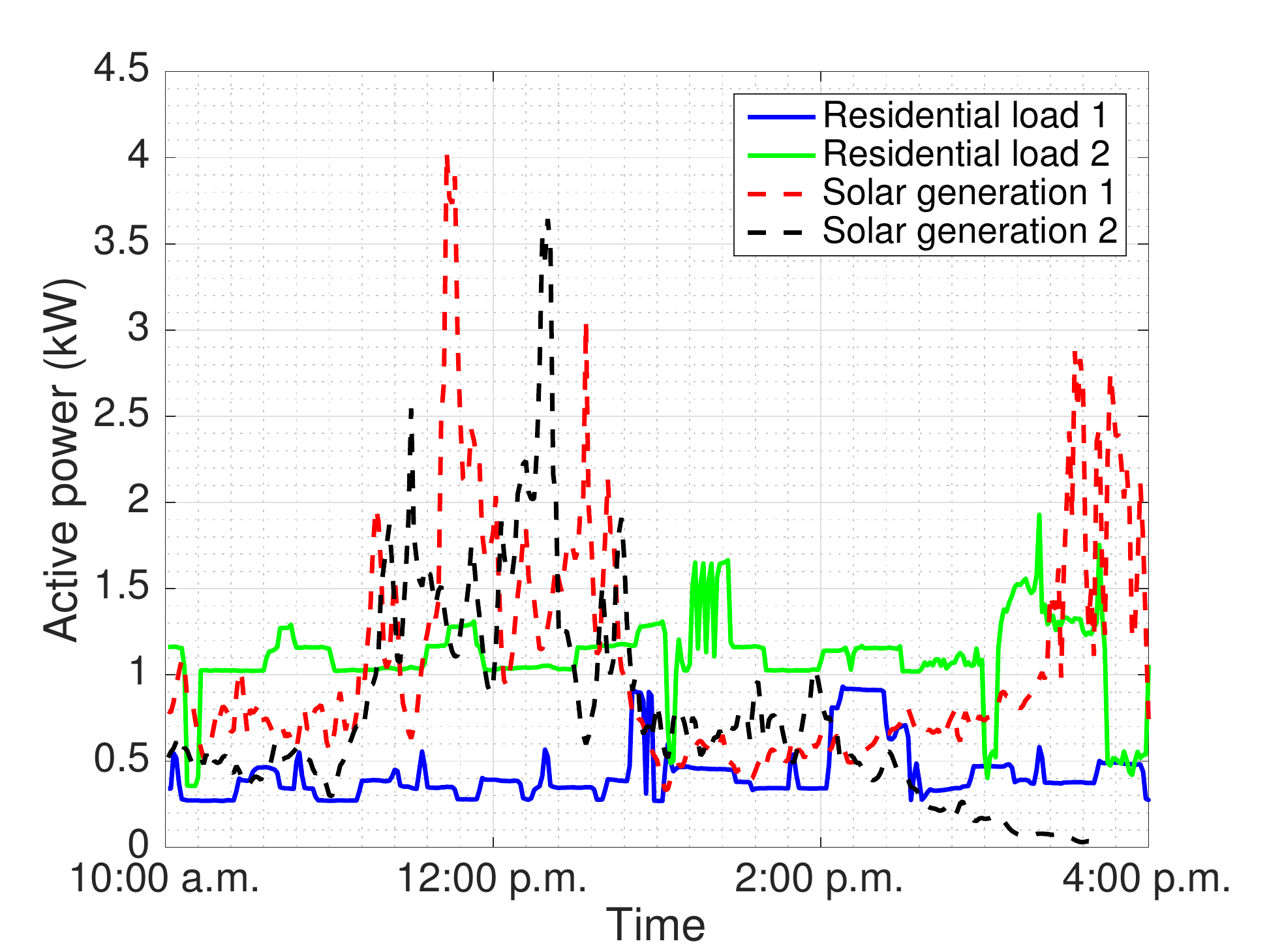}
	\caption{Solar and load data from the Pecan Street project~\cite{pecandata}.}
	\label{fig:pecan}
\end{figure}

To introduce additional specifications, we assume the power injections at non-metered buses $\mathcal{O}$ to remain unchanged between the two time instances. Such an assumption could be justified if buses in $\mathcal{O}$ accommodate conventional loads that are varying slowly compared to the loads in $\mathcal{M}$ or $\mathcal{C}$. Such an assumption can be advocated by Fig.~\ref{fig:pecan} depicting actual data from two typical residential demand profiles and two solar generation curves on a cloudy day~\cite{pecandata}. In this scenario, the grid state  transitions from $\mathbf{v}$ to $\mathbf{v}'$ naturally due to variations in solar generation and elastic demand. 

The stationarity assumption can be alternatively met in an active fashion, where the grid is purposefully perturbed for the sake of inferring loads. This actuation scenario can be realized by commanding the power inverters found in solar panels, energy storage units, or electric vehicles. The inverters can be commanded to promptly and temporarily curtail solar energy, (dis-)charge batteries, or modify their power factor. 

Either naturally or in an active setting, the specifications at metered and conventional buses comply with \eqref{eq:cpf1}--\eqref{eq:cpf2}, while the injection invariability at non-metered buses yields the additional specifications
\begin{subequations}\label{eq:cpf3}
\begin{align} 
p_n(\mathbf{v})&=p_n(\mathbf{v}') & \forall n \in \mathcal{O}\label{eq:cpf3:p}\\
q_n(\mathbf{v})&=q_n(\mathbf{v}') &\forall n \in \mathcal{O}\label{eq:cpf3:q}
\end{align}
\end{subequations}
coupling the two system states. Although the number of unknowns in equations \eqref{eq:cpf1}--\eqref{eq:cpf3} has doubled, the coupling  may allow us to solve the previously under-determined single-time instant problems in \eqref{eq:cpf1} and \eqref{eq:cpf2}. If the specifications are noiseless, then the problem of jointly solving equations \eqref{eq:cpf1}--\eqref{eq:cpf3} is henceforth termed the \emph{coupled power flow (CPF) problem}. 

By simply counting unknowns and equations, a necessary condition for solving \eqref{eq:cpf1}--\eqref{eq:cpf3} is $6M+4C+2O+4\geq 4N+4$, or simply $M\geq O$. A necessary and sufficient condition for local identifiability of the CPF equations is presented  next.

\section{Local Observability Analysis} \label{sec:analysis}
According to the inverse function theorem, a necessary and sufficient condition for locally solving the system of nonlinear equations in \eqref{eq:cpf1}--\eqref{eq:cpf3} is that the related Jacobian matrix $\mathbf{J}(\mathbf{v},\mathbf{v}')$ is invertible. A criterion for characterizing the invertibility of $\mathbf{J}(\mathbf{v},\mathbf{v}')$ is established in this section for the commonly met case of \emph{radial grids}. Some results linking linear algebra to graph theory are reviewed first.

\begin{figure*}[t]
\includegraphics[scale=0.3]{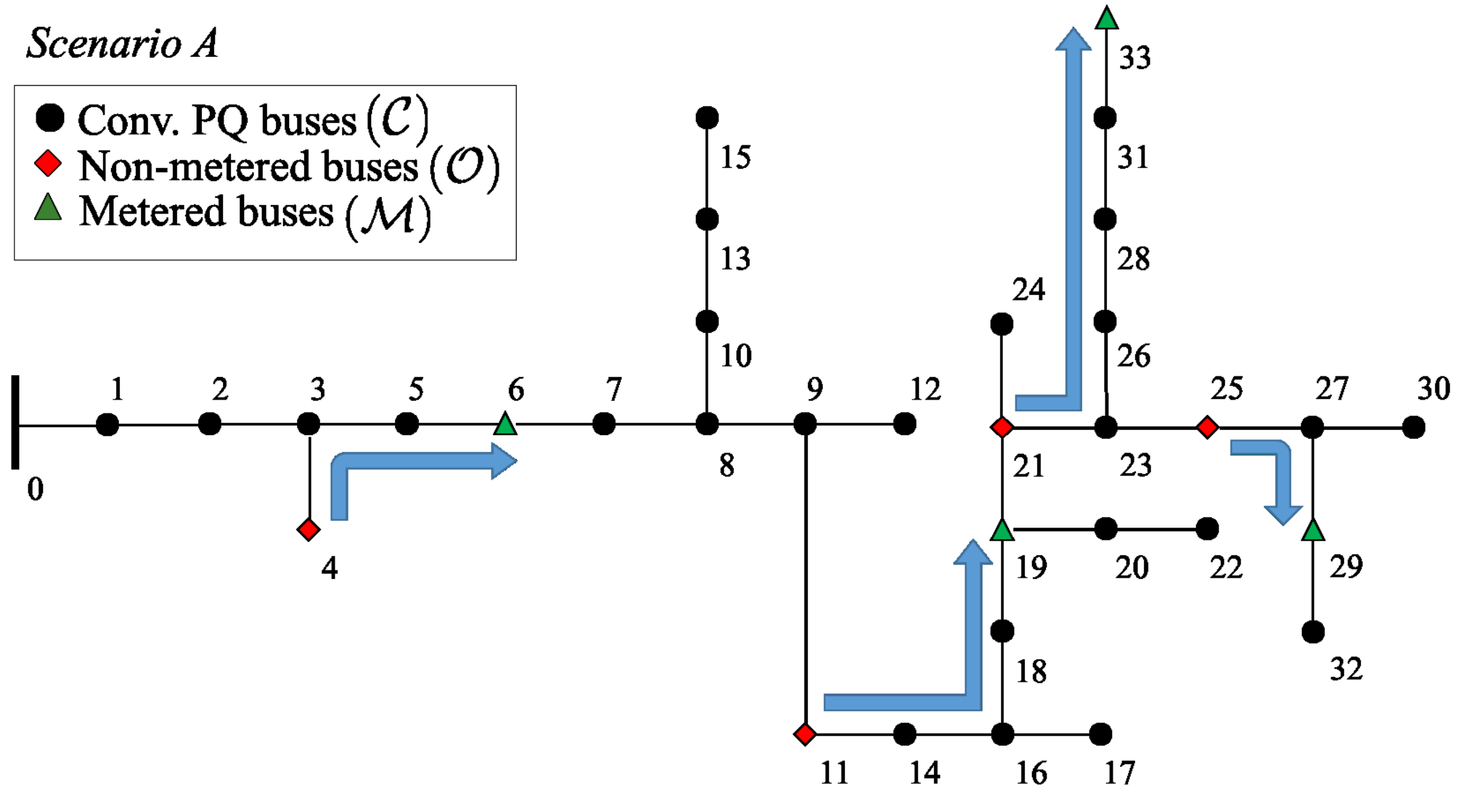}
\hspace*{2em}
\includegraphics[scale=0.3]{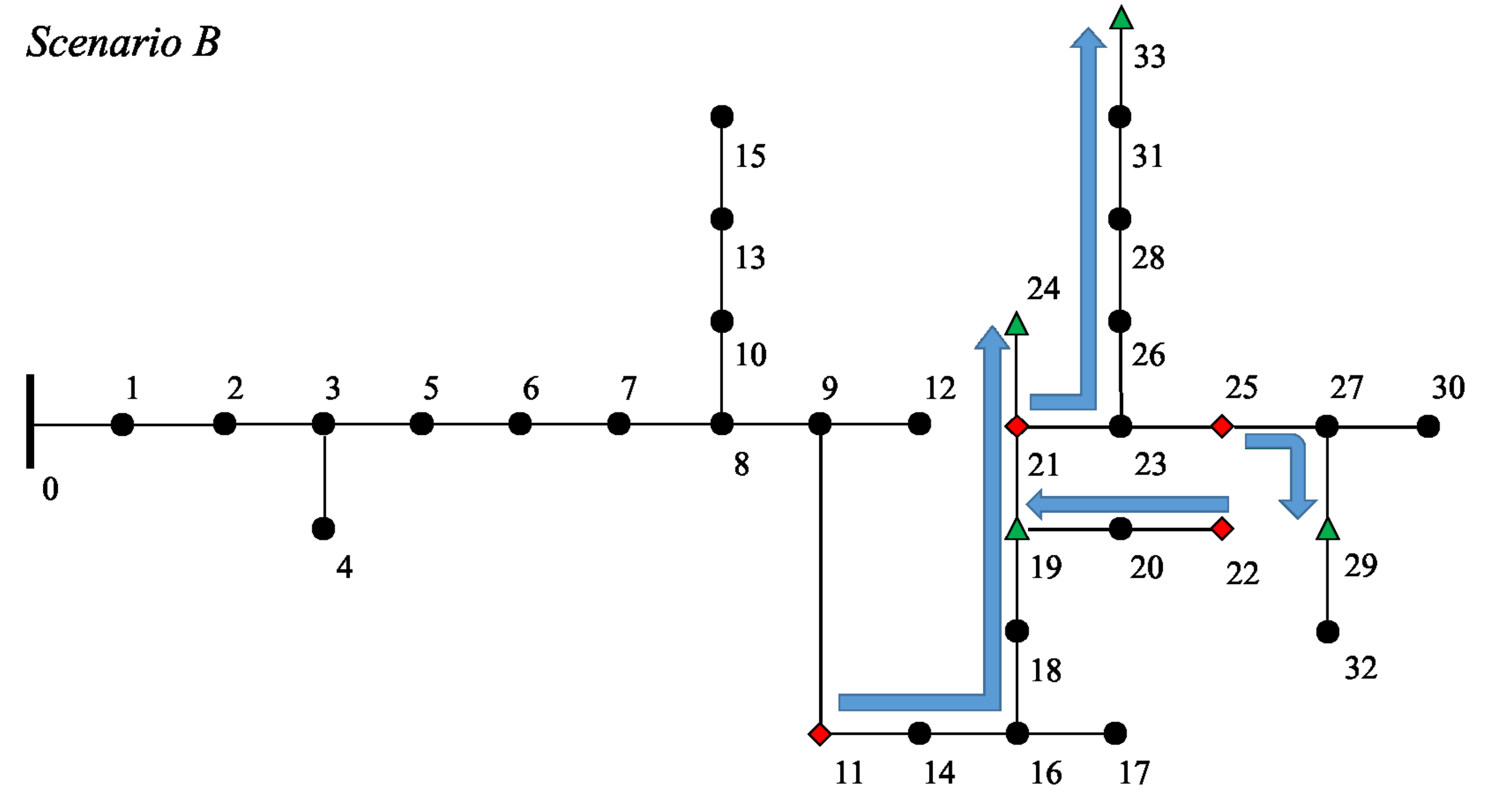}
\caption{Two operating scenarios for the IEEE 34-bus grid~\cite{testfeeder}: the left one satisfies Criterion~\ref{cr:main}; the right one does not.}
\label{fig:34bus}
\end{figure*}

\subsection{Graph Connectivity for Structural Invertibility}\label{subsec:graphs}
A graph $\mathcal{G} = (\mathcal{V}, \mathcal{E})$ is \emph{bipartite} if its vertex set $\mathcal{V}$ can be partitioned into two mutually exclusive and collectively exhaustive subsets $\mathcal{V}_1$ and $\mathcal{V}_2$ such that every edge $e \in \mathcal{E}$ connects a node in $\mathcal{V}_1$ with a node in $\mathcal{V}_2$. Moreover, a \emph{matching} $\mathcal{Z}\subseteq \mathcal{E}$ is a subset of edges in $\mathcal{G}$ so that each node in $\mathcal{V}$ appears in at most one edge in $\mathcal{Z}$. A matching is \emph{perfect} if each node in $\mathcal{V}$ has exactly one edge incident to it. A set of paths in $\mathcal{G}$ is termed \emph{vertex-disjoint} if no two paths share any vertex. 

We are interested in the \emph{generic rank} of a matrix, defined as the maximum possible rank attained if the non-zero entries of this matrix are allowed to take arbitrary values in $\mathbb{R}$. An $N\times N$ matrix is termed generically invertible if its generic rank is $N$. The generic rank of a matrix depends only on its sparsity pattern, that is the locations of its zero entries. 

Let us associate matrix $\mathbf{E}\in\mathbb{R}^{N\times N}$ with a bipartite graph $\mathcal{G}$ having $2N$ nodes. Each column of $\mathbf{E}$ is mapped to a column node, and each row of $\mathbf{E}$ to a row node. An edge runs from the $n$-th column node to the $m$-th row node only if $E_{m,n}\neq 0$. Based on $\mathcal{G}$, the ensuing result holds~\cite{Tutte}.

\begin{lemma}[\cite{Tutte}]\label{le:bipartite}
Matrix $\mathbf{E}$ is generically invertible if and only if the bipartite graph $\mathcal{G}$ has a perfect matching.
\end{lemma}

The result extends to block-partitioned matrices~\cite{VWoude2}. Partition matrix $\mathbf{E}$ as
\begin{equation}\label{eq:pmatrix}
\mathbf{E}=
\begin{bmatrix}
\mathbf{A} & \mathbf{B} \\
\mathbf{C} & \mathbf{D}
\end{bmatrix}
\end{equation}
where $\mathbf{A} \in \mathbb{R}^{M\times M}$, $\mathbf{B} \in \mathbb{R}^{M\times K}$, $\mathbf{C} \in \mathbb{R}^{T\times M}$, and $\mathbf{D} \in \mathbb{R}^{T\times K}$ with $T\geq K$. To study the generic rank of $\mathbf{E}$, construct a directed graph $\mathcal{G}_p=(\mathcal{V},\mathcal{E})$ as follows. First, introduce:
\begin{itemize}
\item a set of state vertices $\mathcal{X}:=\{x_1, \ldots, x_M\}$; 
\item a set of input vertices $\mathcal{U}:=\{u_1, \ldots, u_K\}$; 
\item a set of output vertices $\mathcal{Y}:=\{y_1, \ldots, y_T\}$ 
\end{itemize}
so that $\mathcal{V}=\mathcal{X}\cup \mathcal{U}\cup \mathcal{Y}$. Second, define directed edges as:
\begin{itemize} 
	\item $\mathcal{E}_{xx}=\{(x_{m_2},x_{m_1}):A_{m_1m_2}\neq 0\}$;
	\item $\mathcal{E}_{ux}=\{(u_k,x_m):B_{mk}\neq 0\}$;
	\item $\mathcal{E}_{xy}=\{(x_m,y_t):C_{tm} \neq 0\}$;
	\item $\mathcal{E}_{uy}=\{(u_k,y_t):D_{tk} \neq 0\}$;
\end{itemize}
with $\mathcal{E}=\mathcal{E}_{xx} \cup \mathcal{E}_{ux} \cup \mathcal{E}_{xy} \cup \mathcal{E}_{uy}$.

\begin{lemma}[\cite{VWoude2}]\label{le:struct}
Consider the directed graph $\mathcal{G}_p$ associated with matrix $\mathbf{E}$ constructed as detailed earlier. If $\mathbf{A}$ is invertible, the generic rank of the Schur complement $\mathbf{D}-\mathbf{C}\mathbf{A}^{-1}\mathbf{B}$ is equal to the maximal number of vertex-disjoint paths from $\mathcal{U}$ to $\mathcal{Y}$ in $\mathcal{G}_p$.
\end{lemma}

To characterize the generic rank of $\mathbf{E}$, the rank additivity property for partitioned matrices can be invoked~\cite[p.~25]{HornJohnson}:

\begin{lemma}[Rank additivity]\label{le:rap}
If matrix $\mathbf{A}$ in \eqref{eq:pmatrix} is invertible, it holds that \[\rank(\mathbf{E})= \rank(\mathbf{A})+\rank(\mathbf{D}-\mathbf{C}\mathbf{A}^{-1}\mathbf{B}).\]
\end{lemma}

From Lemmas~\ref{le:struct} and \ref{le:rap}, it readily follows:

\begin{corollary}[\cite{VWoude2}]\label{co:sr}
	If $\mathbf{A}$ in \eqref{eq:pmatrix} is invertible, the generic rank of $\mathbf{E}$ is equal to $M$ plus the maximal number of vertex-disjoint paths from $\mathcal{U}$ to $\mathcal{Y}$ in $\mathcal{G}_p$.
\end{corollary}

\subsection{A Graph-Theoretic Criterion}\label{subsec:criterion}
The next graph-theoretic criterion is critical for the invertibility of the CPF equations.

\begin{criterion}\label{cr:main}
Assume without loss of generality that the graph derived from the tree $\mathcal{G}_d=(\mathcal{N}^+,\mathcal{L})$ by removing the substation bus and its incident edges remains a tree. Nodes are classified into the three mutually exclusive and collectively exhaustive sets of metered $\mathcal{M}$, non-metered $\mathcal{O}$, and conventional nodes $\mathcal{C}$. Assume there exists a vertex-disjoint set of paths connecting every node in $\mathcal{O}$ with a node in $\mathcal{M}$.
\end{criterion}

The main claim on the solvability of the CPF problem is:

\begin{theorem}\label{th:main}
The Jacobian matrix associated with the coupled power flow equations in \eqref{eq:cpf1}--\eqref{eq:cpf3} is generically invertible if and only if Criterion~\ref{cr:main} is satisfied.
\end{theorem}

To better understand the invertibility conditions posed by Theorem~\ref{th:main}, let us study some operational scenarios on the IEEE 34-bus benchmark feeder~\cite{testfeeder}. Consider first the simple scenario where sets $\mathcal{O}$ and $\mathcal{M}$ consist of a single bus each. A unique path can connect these two buses, and there is no other path to be considered. Therefore, if for a single non-metered bus a metered bus is added, the involved Jacobian matrix is invertible in general. Let us study now more complex scenarios with four non-metered buses. Suppose $\mathcal{M}=\{6,19,29,33\}$ and $\mathcal{O}=\{4,11,21,25\}$ shown as Scenario A in Fig.~\ref{fig:34bus}. Since the paths $(4,3,5,6)$, $(11,14,16,18,19)$, $(21,23,26,28,31,33)$ and $(25,27,29)$ do not share any vertex, the related CPF problem is invertible. Consider now the case where $\mathcal{M}=\{19,24,29,33\}$ and $\mathcal{O}=\{11,21,22,25\}$ shown as Scenario B in Fig.~\ref{fig:34bus}. Although bus $25$ can be connected to bus $29$, there is no way to pair $\{11,21,22\}$ to $\{19,24,33\}$ via paths without sharing a common vertex; and hence, this case is not invertible.

Note that for small networks it may be possible to identify the maximum number of vertex-disjoint paths from $\mathcal{O}$ to $\mathcal{M}$ simply by inspection. For larger networks, the original bipartite matching problem can be formulated as a maxflow problem~\cite{Ford}. The latter is typically solved using the Ford-Fulkerson algorithm, whose running time scales linearly with the number of graph nodes and edges~\cite{Ford}.

\subsection{Invertibility of the Jacobian Matrix}\label{subsec:JMI}
Define the mappings of the state $\mathbf{v}$ to the vectors of squared voltage magnitudes and power injections at all buses:
\begin{subequations}\label{eq:map}
	\begin{align}
		\mathbf{m}(\mathbf{v})&:=[|V_0(\mathbf{v})|^2~\ldots~|V_N(\mathbf{v})|^2]^\top \label{eq:map:v}\\
		\mathbf{q}(\mathbf{v})&:=[q_0(\mathbf{v})~\ldots~q_N(\mathbf{v})]^\top\\
		\mathbf{p}(\mathbf{v})&:=[p_0(\mathbf{v})~\ldots~p_N(\mathbf{v})]^\top.\label{eq:map:p}
	\end{align}
\end{subequations}
Nodal power injections are defined as $\mathbf{p}+j\mathbf{q}= \diag(\mathbf{v})\mathbf{i}^*$ with $\mathbf{i}=\mathbf{Y}\mathbf{v}$ being the vector of nodal currents. Eliminating $\mathbf{i}$ yields $\mathbf{p}+j\mathbf{q}= \diag(\mathbf{v})\mathbf{Y}^*\mathbf{v}^*$. Upon differentiating the latter, it can be verified that the Jacobian matrices for the mappings in \eqref{eq:map} can be expressed as~\cite{MLB15}:
\begin{subequations}  \label{eq:Jacobians}
	\begin{align*}
		&\mathbf{J}^m(\mathbf{v}){=}
		\begin{bmatrix}
			2\diag(\mathbf{v}_r) & 2\diag(\mathbf{v}_i)
		\end{bmatrix} 
		\\
		&\mathbf{J}^{q}(\mathbf{v}){=}
		\begin{bmatrix}
			\mathbf{B}\diag({\mathbf{v}_r}){-}\mathbf{G}\diag({\mathbf{v}_i}) &  \mathbf{B}\diag({\mathbf{v}_i}){+}\mathbf{G}\diag({\mathbf{v}_r}) \\
			{+}\diag({\mathbf{Bv}_r}){+}\diag({\mathbf{Gv}_i}) & {-}\diag({\mathbf{Gv}_r}){+}\diag({\mathbf{Bv}_i})
		\end{bmatrix} 
		\\
		&\mathbf{J}^{p}(\mathbf{v}){=}
		\begin{bmatrix}
			{-}\mathbf{G}\diag({\mathbf{v}_r}){-}\mathbf{B}\diag({\mathbf{v}_i}) & {-}\mathbf{G}\diag({\mathbf{v}_i}){+}\mathbf{B}\diag({\mathbf{v}_r}) \\
			{-}\diag({\mathbf{Gv}_r}){+}\diag({\mathbf{Bv}_i}) & {-}\diag({\mathbf{Bv}_r}){-}\diag({\mathbf{Gv}_i})
		\end{bmatrix}.
	\end{align*}
\end{subequations}

The CPF equations can be compactly written as $\mathbf{s(\mathbf{v},\mathbf{v}')}=\mathbf{0}$ with the mapping $\mathbf{s(\mathbf{v},\mathbf{v}')}:\mathbb{R}^{4N+4}\rightarrow \mathbb{R}^{4N+4}$ defined as
\begin{equation} \label{eq:svv}
	\mathbf{s(\mathbf{v},\mathbf{v}')}=\begin{array}{c@{}l}
	\left[\begin{array}{c}
		\mathbf{m}_{\mathcal{M}}(\mathbf{v})-\hat{\mathbf{m}}_{\mathcal{M}}\\
		\mathbf{q}_{\mathcal{C}\cup \mathcal{M}}(\mathbf{v})-\hat{\mathbf{q}}_{\mathcal{C}\cup \mathcal{M}}\\
		\mathbf{p}_{\mathcal{C}\cup \mathcal{M}}(\mathbf{v}) -\hat{\mathbf{p}}_{\mathcal{C}\cup \mathcal{M}}\\
		\mathbf{p}_{\mathcal{O}}(\mathbf{v})-\mathbf{p}_{\mathcal{O}}(\mathbf{v}')\\
		\mathbf{q}_{\mathcal{O}}(\mathbf{v})-\mathbf{q}_{\mathcal{O}}(\mathbf{v}')\\
		\mathbf{m}_{\mathcal{M}}(\mathbf{v}')-\hat{\mathbf{m}}_{\mathcal{M}}'\\
		\mathbf{q}_{\mathcal{C}\cup \mathcal{M}}(\mathbf{v}')-\hat{\mathbf{q}}_{\mathcal{C}\cup \mathcal{M}}'\\
		\mathbf{p}_{\mathcal{C}\cup \mathcal{M}}(\mathbf{v}')-\hat{\mathbf{p}}_{\mathcal{C}\cup \mathcal{M}}'.
	\end{array}\right]
&
\begin{array}[c]{@{}l@{\,}l}
\left. \begin{array}{c} \vphantom{0}  \\ \vphantom{0}
\\ \vphantom{0} \end{array} \right\} & \text{time $t$} \\
\left. \begin{array}{c} \vphantom{0} \\ \vphantom{0}
 \end{array} \right\} & \text{coupled} \\
\left. \begin{array}{c} \vphantom{0}
\\ \vphantom{0} \\ \vphantom{0} \end{array} \right\} & \text{time $t$'} 
\end{array}
\end{array}
\end{equation}
The notation $\mathbf{a}_{\mathcal{B}}$ denotes the subvector of $\mathbf{a}$ indexed by set $\mathcal{B}$. For ease of exposition, the equations in \eqref{eq:cpf3} coupling the two states $(\mathbf{v},\mathbf{v}')$ have been inserted between equations \eqref{eq:cpf1} and \eqref{eq:cpf2} pertaining to single states. Based on the Jacobian matrices for the mappings in \eqref{eq:map}, the Jacobian matrix of $\mathbf{s(\mathbf{v},\mathbf{v}')}$ is
\begin{equation}\label{eq:Jvv}
	\mathbf{J}(\mathbf{v},\mathbf{v}')=
	\left[
	\begin{array}{c:c}
		\mathbf{J}_{\mathcal{M}}^m(\mathbf{v}) & \mathbf{0}\\
		\mathbf{J}_{\mathcal{C} \cup \mathcal{M}}^{q}(\mathbf{v}) & \mathbf{0}\\
		\mathbf{J}_{\mathcal{C} \cup \mathcal{M}}^{p}(\mathbf{v}) & \mathbf{0}\\
		\mathbf{J}_{\mathcal{O}}^{p}(\mathbf{v}) & -\mathbf{J}_{\mathcal{O}}^{p}(\mathbf{v}') \\
		\hdashline
		\mathbf{J}_{\mathcal{O}}^{q}(\mathbf{v}) & -\mathbf{J}_{\mathcal{O}}^{q}(\mathbf{v}') \\
		\mathbf{0} & \mathbf{J}_{\mathcal{M}}^m(\mathbf{v}') \\
		\mathbf{0} & \mathbf{J}_{\mathcal{C} \cup \mathcal{M}}^{q}(\mathbf{v}')\\
		\mathbf{0} & \mathbf{J}_{\mathcal{C} \cup \mathcal{M}}^{p}(\mathbf{v}')\\
	\end{array}
	\right]
\end{equation}
that can be partitioned as
\begin{equation}\label{eq:Jvvp}
	\mathbf{J}(\mathbf{v},\mathbf{v}')=
	\left[
	\begin{array}{c:c}
		\mathbf{J}_A(\mathbf{v}) & \mathbf{J}_B(\mathbf{v}') \\
		\hdashline
		\mathbf{J}_C(\mathbf{v}) & \mathbf{J}_D(\mathbf{v}')
	\end{array}
	\right].
\end{equation}

Since studying the rank of $\mathbf{J}(\mathbf{v},\mathbf{v}')$ for any $(\mathbf{v},\mathbf{v}')$ is challenging, we resort to its generic rank instead. To this end, Lemma~\ref{le:rap} is used: We first show the generic invertibility of $\mathbf{J}_A(\mathbf{v})$ and proceed with the generic invertibility of the related Schur complement $\mathbf{J}_{D}(\mathbf{v}')-\mathbf{J}_{C}(\mathbf{v})\mathbf{J}_{A}^{-1}(\mathbf{v})\mathbf{J}_{B}(\mathbf{v}')$. Both results are proved in the Appendix.

\begin{lemma}[Invertibility of $\mathbf{J}_A(\mathbf{v})$]\label{le:JA}
Under Criterion~\ref{cr:main}, the partial Jacobian matrix $\mathbf{J}_A(\mathbf{v})$ is generically invertible.
\end{lemma}

\begin{remark}
Note that since $\mathbf{J}_A(\mathbf{v})$ is generically invertible, the  single-state problem can be solved under the same criterion if a type $P$ or $Q$ specification is available for all buses in $\mathcal{O}$.
\end{remark}

Criterion~\ref{cr:main} further guarantees the invertibility of the Schur complement of $\mathbf{J}(\mathbf{v},\mathbf{v}')$ with respect to $\mathbf{J}_{A}$:

\begin{lemma}[Invertibility of Schur complement]\label{le:Schur}
Under Criterion~\ref{cr:main}, the Schur complement $\mathbf{J}_{D}(\mathbf{v}')-\mathbf{J}_{C}(\mathbf{v})\mathbf{J}_{A}^{-1}(\mathbf{v})\mathbf{J}_{B}(\mathbf{v}')$ is generically invertible.
\end{lemma}

Theorem~\ref{th:main} follows readily from Lemmas~\ref{le:JA} and \ref{le:Schur}. Even though Theorem~\ref{th:main} characterizes only the generic rank of $\mathbf{J(v,v')}$ to be full, the matrix was actually invertible for different state pairs tested in Section~\ref{sec:sim} conditioned on Criterion~\ref{cr:main}.

\section{Solving the Coupled Power Flow Problem}\label{sec:solvers}
Having characterized the local observability for CPF, solvers using exact and noisy specifications are presented next.

\subsection{Noiseless Specifications}\label{sec:cpf}
If grid specifications are exact and the power injections in $\mathcal{O}$ remain unaltered between $\mathbf{v}$ and $\mathbf{v}'$, the CPF task boils down to solving \eqref{eq:cpf1}--\eqref{eq:cpf3}. The classic PF problem can be solved using the Newton-Raphson (NR) method, its decoupled variants, Gauss-Seidel, or Jacobi iterations~\cite{ExpConCanBook}. Moreover, the forward-backward sweep algorithm has been widely used in radial distribution grids \cite{ExpConCanBook}. Nevertheless, the latter cannot be used for CPF due to the coupling equations. The NR scheme while applicable to CPF is known to diverge if not properly initialized and a good initial guess may not be available with high renewable penetration. For this reason, we address our CPF formulation by adopting the SDP-based PF solver developed in \cite{MLB15} and \cite{MALB16}, and briefly reviewed next.

The key idea of the SDP-based solver is that the PF specifications in \eqref{eq:PF0}--\eqref{eq:PF} are quadratic functions of $\mathbf{v}$ and therefore can be expressed as $\mathbf{v}^H\mathbf{M}_k\mathbf{v}=\hat{s}_k$ for a specific Hermitian matrix $\mathbf{M}_k$ for all $k$. By introducing the matrix variable $\mathbf{V}=\mathbf{vv}^H$, the specifications can be equivalently written as $\trace(\mathbf{M}_k\mathbf{V})=\hat{s}_k$ and the PF task can be posed as the feasibility problem:
\begin{subequations} \label{eq:sdpf1}
	\begin{align}
		\find~&~(\mathbf{V},\mathbf{v}) \label{eq:sdpc2}\\
		\textrm{s.t.}~&~\trace(\mathbf{M}_k\mathbf{V})=\hat{s}_k,\quad k = 1, \ldots, 2N+2 \label{eq:sdpc3}\\
		~&~\mathbf{V}=\mathbf{vv}^H. \label{eq:sdpc4}
	\end{align}
\end{subequations}
Although the constraints in \eqref{eq:sdpc3} are linear with respect to $\mathbf{V}$, the constraint in \eqref{eq:sdpc4} is non-convex. To see this, note that $\mathbf{V}=\mathbf{vv}^H$ is equivalent to $\mathbf{V} \succeq \mathbf{0}$ and $\rank(\mathbf{V})=1$. Problem \eqref{eq:sdpf1} can be relaxed by dropping the non-convex rank constraint and replacing $\mathbf{V}=\mathbf{vv}^H$ with $\mathbf{V} \succeq \mathbf{0}$. To locate feasible points of rank-one, the feasibility problem is turned into the SDP minimization task~\cite{MLB15}: 
\begin{subequations} \label{eq:SDP1}
	\begin{align}
		\min_{\mathbf{V} \succeq \mathbf{0}} ~&~ \trace(\mathbf{M}\mathbf{V})\\
		\textrm{s.t.} ~&~\trace(\mathbf{M}_k\mathbf{V})=\hat{s}_k,\quad k = 1, \ldots, 2N+2.
	\end{align}
\end{subequations}
The matrix $\mathbf{M}$ is judiciously selected to provide rank-one minimizers if the power system operates close to the flat voltage profile. 
If problem \eqref{eq:SDP1} yields a rank-one minimizer with eigenvalue decomposition $\mathbf{V}^*=\lambda\mathbf{uu}^H$, a PF solution is recovered as $\mathbf{v}=\sqrt{\lambda}\mathbf{u}$. 

Building on \eqref{eq:SDP1}, the CPF task can be posed as:
\begin{subequations} \label{eq:CPF}
	\begin{align}
		\min_{\mathbf{V}, \mathbf{V}'  \succeq \mathbf{0}} ~&~ \trace(\mathbf{M}\mathbf{V}) + \trace(\mathbf{M}\mathbf{V}') \label{eq:CPFa}\\
		\textrm{s.t.} ~&~ \trace(\mathbf{M}_{vn}\mathbf{V})=|v_{n}|^2, &\forall n \in \mathcal{S}\cup\mathcal{M} \label{eq:CPFb}\\
		~&~ \trace(\mathbf{M}_{pn}\mathbf{V})=p_{n},&\forall n \in \mathcal{C}\cup\mathcal{M} \label{eq:CPFc}\\
		~&~ \trace(\mathbf{M}_{qn}\mathbf{V})=q_{n}, &\forall n \in \mathcal{C}\cup\mathcal{M} \label{eq:CPFd}\\
		~&~ \trace(\mathbf{M}_{vn}\mathbf{V}')=|v'_{n}|^2, &\forall n \in \mathcal{S}\cup\mathcal{M} \label{eq:CPFe}\\
		~&~\trace(\mathbf{M}_{pn}\mathbf{V}')=p'_{n},&\forall n \in \mathcal{C}\cup\mathcal{M} \label{eq:CPFf}\\
		~&~\trace(\mathbf{M}_{qn}\mathbf{V}')=q'_{n},&\forall n \in  \mathcal{C}\cup\mathcal{M} \label{eq:CPFg}\\
		~&~  \trace(\mathbf{M}_{pn}\mathbf{V})=\trace(\mathbf{M}_{pn}\mathbf{V}'), &\forall n \in \mathcal{O} \label{eq:CPFh}\\
		~&~ \trace(\mathbf{M}_{qn}\mathbf{V})=\trace(\mathbf{M}_{qn}\mathbf{V}'),&\forall n \in \mathcal{O} \label{eq:CPFi}
	\end{align}
\end{subequations}
where $\mathbf{M}_{vn}:= \mathbf{e}_n\mathbf{e}^H_n$, $\mathbf{M}_{pn} := \frac{1}{2}(\mathbf{Y}^H\mathbf{e}_n\mathbf{e}_n^H + \mathbf{e}_n\mathbf{e}_n^H\mathbf{Y})$, $\mathbf{M}_{qn} := \frac{1}{2j}(\mathbf{Y}^H\mathbf{e}_n\mathbf{e}_n^H - \mathbf{e}_n\mathbf{e}_n^H\mathbf{Y})$, and $\mathbf{e}_n$ is the $n$-th the canonical vector in $\mathbb{R}^N$. Constraints \eqref{eq:CPFb}--\eqref{eq:CPFd} correspond to the specifications related to state $\mathbf{v}$; constraints \eqref{eq:CPFe}--\eqref{eq:CPFg} to state $\mathbf{v}'$; and \eqref{eq:CPFh}--\eqref{eq:CPFi} capture the coupling specifications. Although Section~\ref{sec:analysis} considered exactly $4(N+1)$ specifications, adding extra specifications as constraints in \eqref{eq:CPF} could only increase the chances of recovering the exact states $(\mathbf{v},\mathbf{v}')$.

\subsection{Solving CPF under Noisy Measurements}\label{subsec:cpsse}
In practice, the coupled power flow problem has to be solved under noisy scenarios. If voltage magnitudes and power injection specifications in $\mathcal{M}$ and $\mathcal{C}$ come from smart meter readings, they involve measurement and possibly modeling noise, while communication delays could also be captured by imprecise specifications in \eqref{eq:cpf1}--\eqref{eq:cpf2}. Moreover, the coupling equations in \eqref{eq:cpf3} will not be exact, due to small perturbations in the unknown conventional load of $\mathcal{O}$.

To address noisy smart meter data, we put forth a PSSE task coupling two successive grid states, henceforth termed \emph{coupled power system state estimation (CPSSE)}. In the presence of noisy measurements, the convex problem in \eqref{eq:CPF} may not yield a good approximate solution, even if the number of specifications are larger than $4(N+1)$. Given that the Gauss-Newton method, conventionally employed in PSSE, exhibits similar behavior to NR, we extend the penalized SDP-based PSSE solver of~\cite{MALB16} to our CPSSE setting as outlined next:
\begin{subequations} \label{eq:CPSSE}
\begin{align}
\min_{\mathbf{V},\mathbf{V}',\boldsymbol{\epsilon}} &~\sum_{\ell=1}^{2L+2O} f_{\ell}(\epsilon_{\ell}) + \alpha[\trace(\mathbf{M}\mathbf{V}) + \trace(\mathbf{M}\mathbf{V}')]\label{eq:CPSSE:cost}\\
\textrm{s.t.} ~&\trace(\mathbf{M}_\ell\mathbf{V})+\epsilon_{\ell}=\hat{s}_\ell, \quad\ell=1,\ldots,L\label{eq:CPSSE:1}\\
~&\trace(\mathbf{M}_\ell\mathbf{V}')+\epsilon_{\ell}=\hat{s}_\ell', \quad\ell=L+1,\ldots,2L\label{eq:CPSSE:2}\\
~&\trace(\mathbf{M}_\ell(\mathbf{V}'-\mathbf{V}))=\epsilon_{\ell}, ~\ell=2L+1,\ldots,2L+2O\label{eq:CPSSE:O}\\
~& \mathbf{V} \succeq \mathbf{0}, \mathbf{V}' \succeq \mathbf{0}.
\end{align}
\end{subequations}
Measurements $\{\hat{s}_\ell\}_{\ell=1}^L$ relate to state $\mathbf{v}$ in \eqref{eq:CPSSE:1}; $\{\hat{s}'_\ell\}_{\ell=1}^L$ relate to $\mathbf{v}'$ in \eqref{eq:CPSSE:2}; and \eqref{eq:CPSSE:O} couples the two states. 

The optimization variables $\{\epsilon_{\ell}\}_{\ell=1}^{2L+2O}$ collected in vector $\boldsymbol{\epsilon}$ capture the residuals between the actual grid quantities and the metered or coupled data. Residual $\epsilon_\ell$ appears in the first summand of \eqref{eq:CPSSE:cost} as the argument of the data-fitting function $f_\ell(\epsilon_\ell)$. This function can be either a weighted squared or absolute value, i.e., $f_{\ell}(\epsilon_\ell)= \epsilon_\ell^2/\sigma_\ell^2$ or $|\epsilon_\ell|/\sigma_\ell$ with different $\sigma_\ell$'s depending on the uncertainty of the $\ell$-th datum.

The second summand in \eqref{eq:CPSSE:cost} corresponds to a regularizer promoting rank-one minimizers with the tuning parameter $\alpha>0$. For $\alpha=0$, the CPSSE cost involves only the data-fitting term. For increasing $\alpha$, more emphasis is placed on the cost components $\trace(\mathbf{M}\mathbf{V}) + \trace(\mathbf{M}\mathbf{V}')$~\cite{MALB16}. If the minimizers $\mathbf{V}^*$ and $(\mathbf{V}')^*$ of \eqref{eq:CPSSE} are not rank-one, the heuristic for finding voltage angles proposed in~\cite{MALB16} is followed. Knowing that non-metered buses host loads, the additional constraints $\trace(\mathbf{M}_{pn}\mathbf{V})<0$ and $\trace(\mathbf{M}_{pn}\mathbf{V}')<0$ for all $n\in \mathcal{O}$ have been appended in \eqref{eq:CPSSE} to strengthen the SDP relaxation.

\section{Numerical Tests}\label{sec:sim}
The solvability criterion for the CPF equations along with the CPF and CPSSE solvers were numerically tested using the IEEE 34-bus feeder. The original multi-phase grid was converted into an equivalent single-phase one using the procedure described in~\cite{GLTL12}. Additionally, at each zero-injection bus, a load equal to the load of its parent node was inserted. The tests were run on a 2.7 GHz Intel Core i5 with 8GB RAM using the SDPT3 solver in YALMIP and MATLAB~\cite{YALMIP}, \cite{SDPT3}.

\begin{figure*}[t]
	\centering
	\includegraphics[scale=0.38]{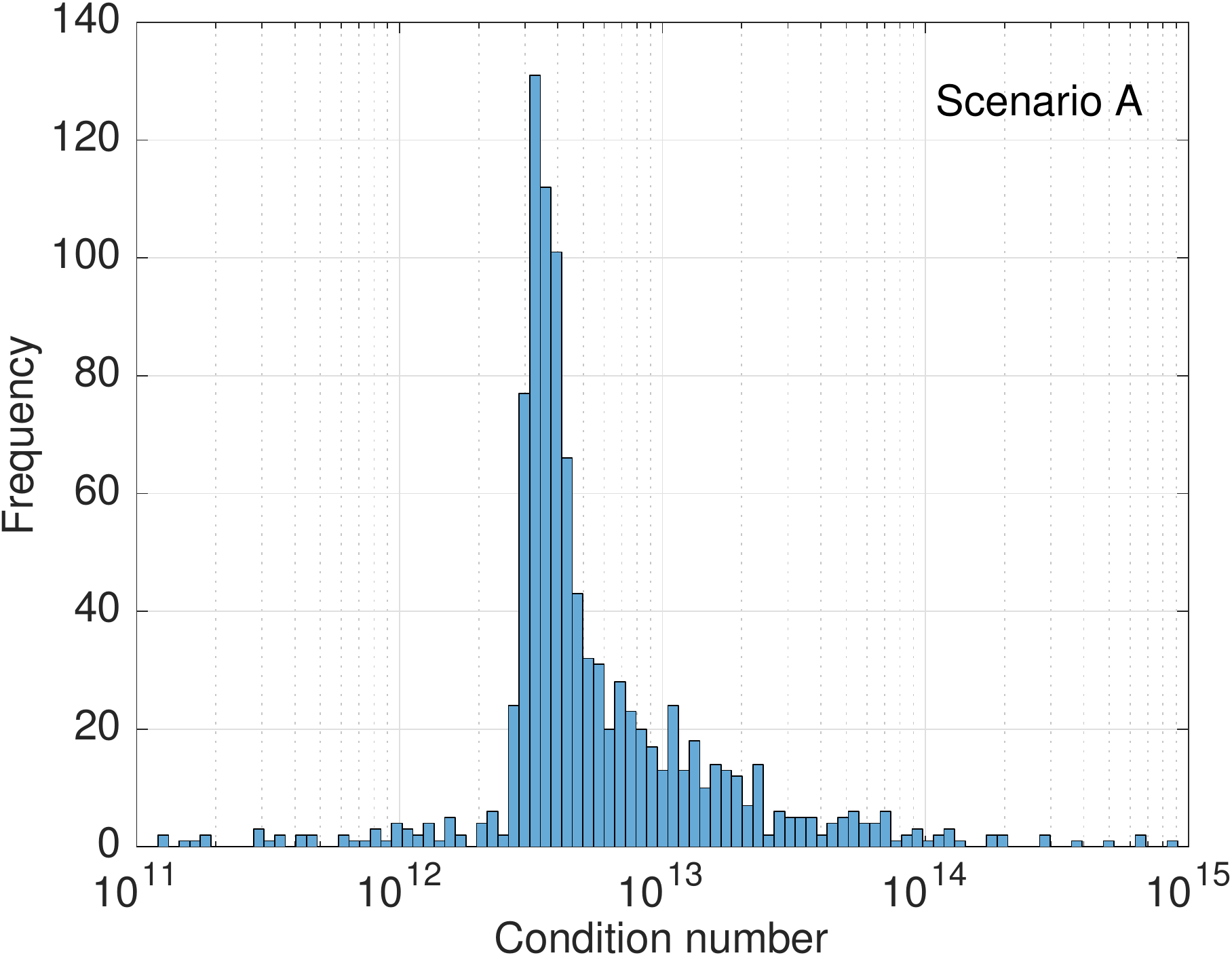}
	\hspace*{4em}
	\includegraphics[scale=0.38]{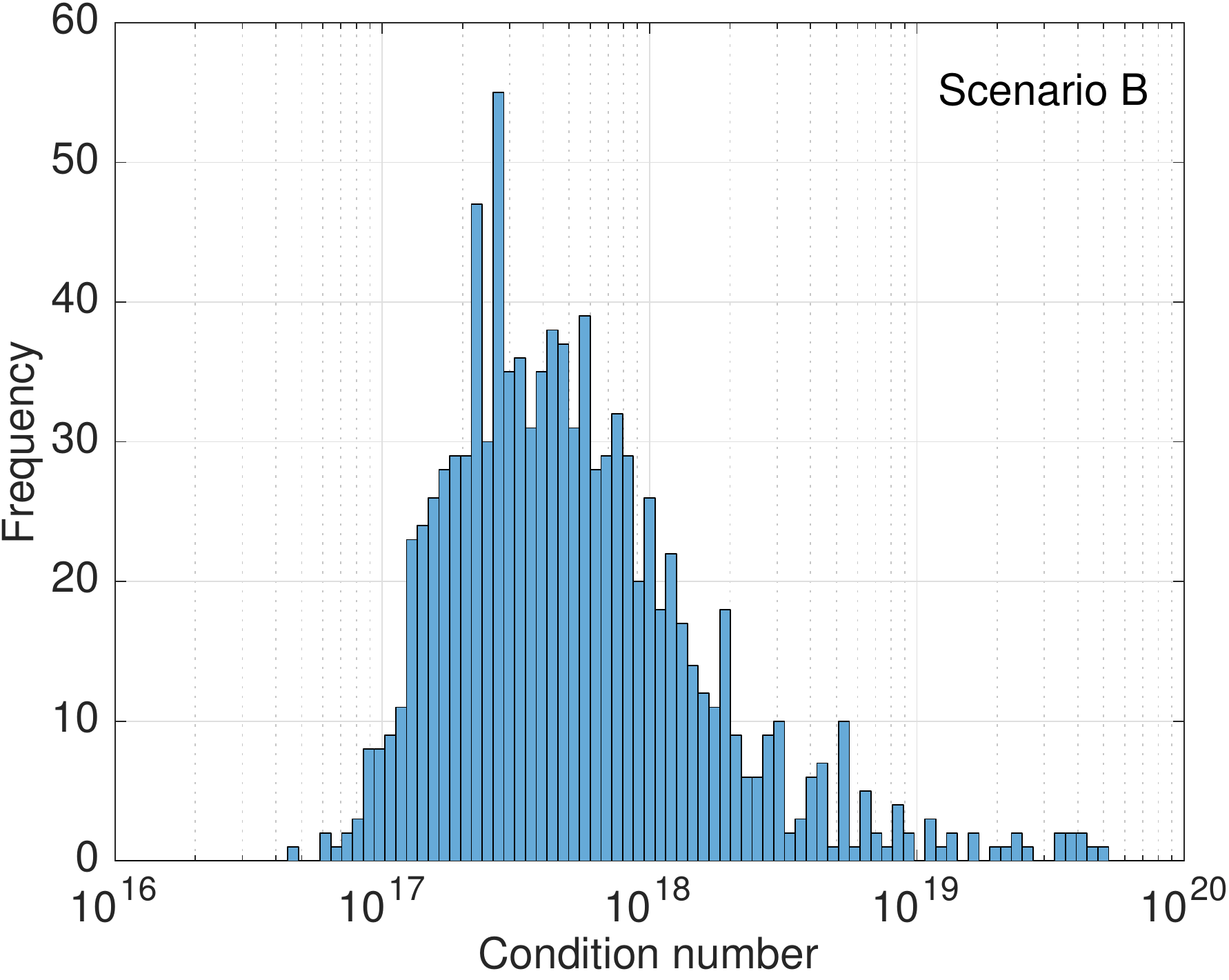}
	\caption{Histograms for the condition numbers of the Jacobian matrix $\mathbf{J}(\mathbf{v},\mathbf{v}')$ for the IEEE 34-bus grid shown in Fig.~\ref{fig:34bus}.}
	\label{fig:Hist}
\end{figure*} 

We first tested the validity of Criterion~\ref{cr:main}. To simulate solar generation, a photovoltaic with capacity 4 times the related load was added on all metered buses and on conventional PQ buses $\{2,9,13,16,24,28\}$. These same buses were assumed controllable. To show that Criterion~\ref{cr:main} characterizes CPF solvability, matrix $\mathbf{J}(\mathbf{v},\mathbf{v}')$ and its condition number was evaluated for 1,000 random states $(\mathbf{v},\mathbf{v}')$. For grid state $\mathbf{v}$, solar generators produced their maximum active power at 0.9 lagging power factor. For state $\mathbf{v}'$, solar generation was uniformly drawn within its capacity range while fixing the power factor to 0.9 leading or lagging. For each realization of $(\mathbf{v},\mathbf{v}')$, Fig.~\ref{fig:Hist} shows the histograms of the condition numbers related to the scenarios of Fig.~\ref{fig:34bus}. As evidenced by the plots, the Jacobian for Scenario A exhibits a condition number much smaller than Scenario B in accordance with Criterion~\ref{cr:main}.

\begin{figure*}[t]
	\centering
	\includegraphics[scale=0.38]{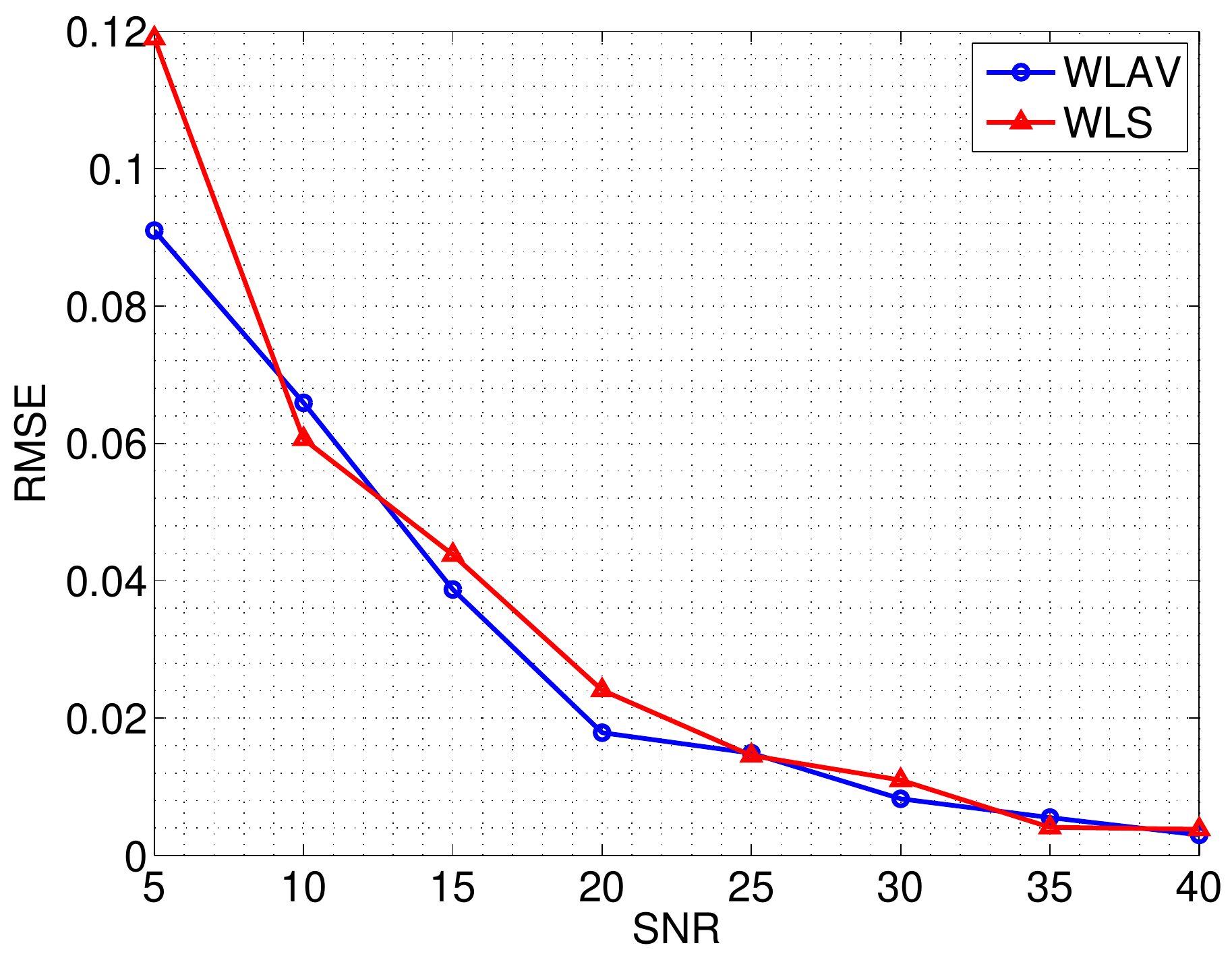}
	\hspace*{4em}
	\includegraphics[scale=0.38]{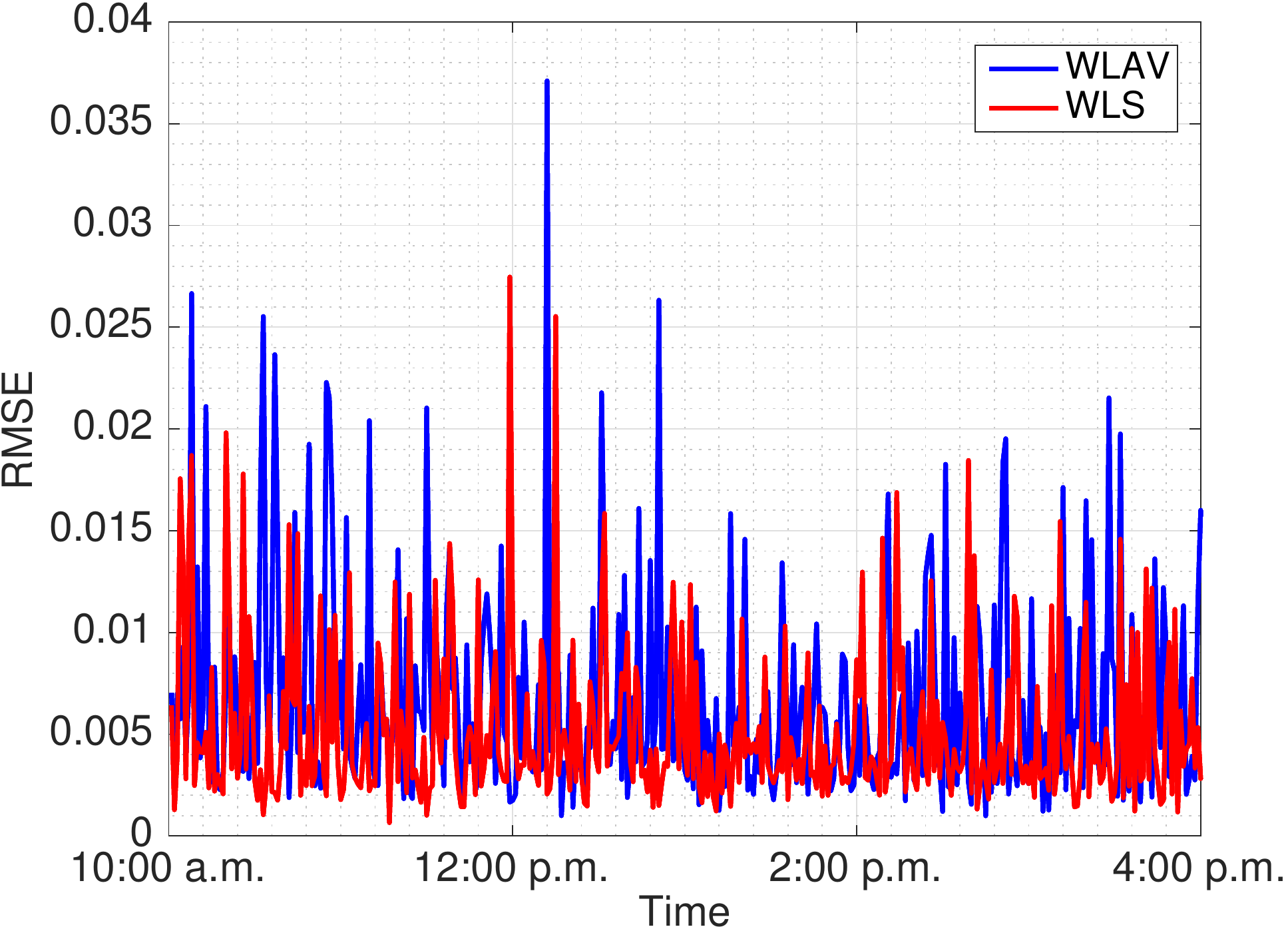}
	\caption{CPSSE using actual load and solar generation~\cite{pecandata}. RMSE for $(\mathbf{v},\mathbf{v}')$ for \eqref{eq:CPSSE} (left) and over a day (right).}
	\vspace*{-1em}
	\label{fig:SDPSE}
\end{figure*} 

Assuming noiseless specifications, the CPF solver of \eqref{eq:CPF} was subsequently tested for Scenario A. The matrix $\mathbf{M}$ in \eqref{eq:CPFa} was set to $-\mathbf{B}$. Solar generators produced again their maximum active power $p_n$ at 0.9 lagging power factor for state $\mathbf{v}$. For state $\mathbf{v}'$, solar generation was perturbed by $\beta z_n p_n$, where $\beta$ ranged between $(0,1]$ in increments of $0.1$ and ${z}_n$ was uniformly drawn in $[0,0.1]$. The power factor was fixed to 0.9 lagging. For each value of $\beta$, 50 realizations of $\mathbf{v}'$ were constructed. The probability of success of the CPF solver is shown in Fig.~\ref{fig:ProbSDP}. As $\beta$ increases, the state $\mathbf{v}'$ moves away from $\mathbf{v}$ and the probability of success increases. When states $\mathbf{v}$ and $\mathbf{v}'$ are relatively close, matrix $\mathbf{J}(\mathbf{v},\mathbf{v}')$ becomes numerically singular and the solver fails as expected.

Given noisy specifications, the CPSSE task was next tackled using actual data and the SDP-based solver of \eqref{eq:CPSSE}. Load and solar generation on the IEEE 34-bus grid was taken from the Pecan Street dataset between 10:00 a.m. and 4:00 p.m. on January 1, 2013~\cite{pecandata}. Solar data from 10 homes were added on metered buses $\{2,6,9,13,16,19,24,28,29,33\}$. Consumption data from 23 different homes were added as spot loads on the remaining buses except for the substation. Data were normalized and scaled by the static peak load of the IEEE-34 bus system. Due to lack of reactive data, power factors were simulated to range between 0.95 leading to 0.95 lagging. 

The left panel of Fig.~\ref{fig:SDPSE} presents the root-mean square error (RMSE) as a function of the measurement signal-to-noise ratio (SNR). The RMSE was averaged over 20 instances of zero-mean Gaussian measurement noise for each SNR value. The right panel of Fig.~\ref{fig:SDPSE} shows the accuracy of the $(\mathbf{v},\mathbf{v}')$ estimates over the day for the weighted least-squares (WLS) and the weighted least absolute value (WLAV) costs. In this test, the measurement noise was simulated as Gaussian with zero mean and standard deviation of $0.01$ for voltages and $0.015$ for power injections. The coupling equations in \eqref{eq:CPSSE} were weighted by $0.035$ and $\alpha$ was set to 2 in all tests.

\section{Conclusions}\label{sec:conclusions}
Smart meter data together with the actuation capabilities of power inverters were leveraged to recover the state of distribution grids. To compensate for missing specifications on non-metered buses, we exploited the relative stationarity of conventional loads and coupled grid states across time. The solvability of the related nonlinear equations was characterized via an intuitive and easily verifiable graph-theoretic criterion. We further put forth solvers relying on SDP relaxations for the noiseless and noisy variants of the joint power flow problem. Inferring the grid state through coupling was shown to be successful using actual residential load and solar generation data on a benchmark feeder. Extending the approach to multiple time instances and unbalanced multi-phase conditions, involving approximate grid models, and designing optimal system actuation schemes constitute our current research efforts.


\appendix\label{sec:appendix}
\begin{IEEEproof}[Proof of Lemma~\ref{le:JA}]
Rather than checking the rank of $\mathbf{J}_A(\mathbf{v})$ for all $\mathbf{v}$, matrix $\mathbf{J}_A(\mathbf{v})$ is evaluated at the flat voltage profile $\mathbf{v}_{\textrm{fl}}:= [\mathbf{1}_{N+1}^{\top} ~\mathbf{0}_{N+1}^{\top}]^{\top}$. This is without loss of generality, because the generic rank of a matrix depends only on its sparsity pattern.
	
Using \eqref{eq:Jacobians}, \eqref{eq:Jvv}, and the fact that $\mathbf{G1}=\mathbf{B1}=\mathbf{0}$, matrix $\mathbf{J}_A(\mathbf{v}_{\text{fl}})$ can be written as
	\begin{equation}\label{eq:JA}
	\mathbf{J}_A(\mathbf{v}_{\textrm{fl}})=
	\begin{bmatrix}
	2\mathbf{I}_{\mathcal{S} \cup \mathcal{M},\mathcal{N}^+} & \mathbf{0}_{\mathcal{S} \cup \mathcal{M},\mathcal{N}} \\
	\mathbf{B}_{\mathcal{C} \cup \mathcal{M},\mathcal{N}^+} & \mathbf{G}_{\mathcal{C} \cup \mathcal{M}, \mathcal{N}} \\
	-\mathbf{G}_{\mathcal{N},\mathcal{N}^+} &\mathbf{B}_{\mathcal{N},\mathcal{N}}
	\end{bmatrix}.
	\end{equation}
Since $\mathcal{N}^+=\mathcal{S}\cup \mathcal{M}\cup \mathcal{C}\cup \mathcal{O}$, matrix $\mathbf{J}_A(\mathbf{v}_{\textrm{fl}})$ is partitioned as:
\begin{equation*}
	\mathbf{J}_A(\mathbf{v}_{\textrm{fl}})=
	\left[
	\begin{array}{c:cc}
	2\mathbf{I}_{\mathcal{S} \cup \mathcal{M},\mathcal{S} \cup \mathcal{M}} & \mathbf{0}_{\mathcal{S} \cup \mathcal{M}, \mathcal{C} \cup \mathcal{O}} & 	\mathbf{0}_{\mathcal{S} \cup \mathcal{M},\mathcal{N}} \\
	\hdashline
	\mathbf{B}_{\mathcal{C} \cup \mathcal{M},\mathcal{S} \cup \mathcal{M}} &\mathbf{B}_{\mathcal{C} \cup \mathcal{M},\mathcal{C} \cup \mathcal{O}} & \mathbf{G}_{\mathcal{C} \cup \mathcal{M},\mathcal{N}}\\
	-\mathbf{G}_{\mathcal{N}, \mathcal{S} \cup \mathcal{M}} &
	-\mathbf{G}_{\mathcal{N}, \mathcal{C} \cup \mathcal{O}} &
	\mathbf{B}_{\mathcal{N},\mathcal{N}} 
	\end{array}
	\right].
\end{equation*}

	\begin{figure}[t]
	\centering
	\includegraphics[scale=0.38]{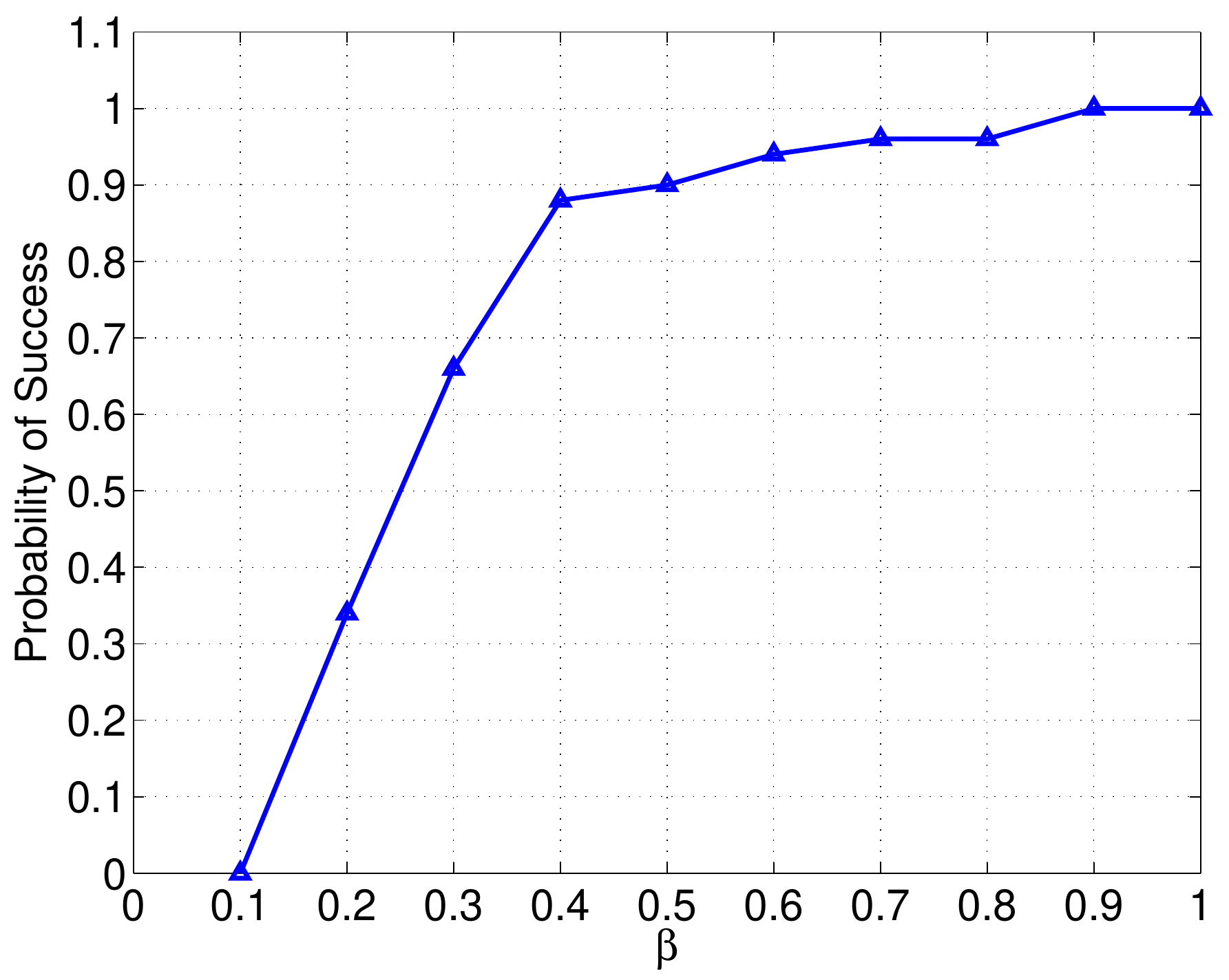}
	\caption{Probability of successful recovery of $(\mathbf{v},\mathbf{v}')$ from \eqref{eq:CPF}.}
	\label{fig:ProbSDP}
\end{figure}

Because the top left block of $\mathbf{J}_A(\mathbf{v}_{\textrm{fl}})$ is a scaled identity matrix and its top right block is zero, Lemma~\ref{le:rap} asserts that $\mathbf{J}_A(\mathbf{v}_{\textrm{fl}})$ is full rank if and only if its bottom right block
	\begin{equation*}
	\bar{\mathbf{J}}_{A}:=
	\begin{bmatrix}
	\mathbf{B}_{\mathcal{C} \cup \mathcal{M},\mathcal{C} \cup \mathcal{O}} & \mathbf{G}_{\mathcal{C} \cup \mathcal{M},\mathcal{N}}\\
	-\mathbf{G}_{\mathcal{N},\mathcal{C} \cup \mathcal{O}} &\mathbf{B}_{\mathcal{N},\mathcal{N}}
	\end{bmatrix}
	\end{equation*}
is full rank. Observe that $\mathbf{B}_{\mathcal{N},\mathcal{N}}$ is a symmetric submatrix of the weighted Laplacian matrix $\mathbf{B}$ and has been shown to be invertible~\cite[Lemma~1]{MLB15}. Using the rank additivity property of Lemma~\ref{le:rap} once again, matrix $\bar{\mathbf{J}}_{A}$ is invertible if and only if the ensuing Schur complement is full rank:
\begin{equation}\label{eq:JAt}
	\tilde{\mathbf{J}}_{A}:= \mathbf{B}_{\mathcal{C} \cup \mathcal{M}, \mathcal{C} \cup \mathcal{O}} + \mathbf{G}_{\mathcal{C} \cup \mathcal{M},\mathcal{N}} \mathbf{B}_{\mathcal{N},\mathcal{N}}^{-1} \mathbf{G}_{\mathcal{N}, \mathcal{C} \cup \mathcal{O}}.
\end{equation}  

To determine the rank of $\tilde{\mathbf{J}}_{A}$, let us express $\mathbf{B}$ and $\mathbf{G}$ in terms of the grid parameters. Let $\mathbf{b}$ and $\mathbf{g}$ be  the $N$--dimensional vectors of distribution line susceptances and conductances, respectively. Moreover, let $\mathbf{A}$ be the reduced bus-branch incidence matrix~\cite{SGC15}. Matrix $\mathbf{A}$ is square and invertible for a radial network. Using these definitions, the susceptance matrix is $\mathbf{B}=\mathbf{A}^{\top}\diag(\mathbf{b}) \mathbf{A}$ and the conductance matrix is $\mathbf{G}=\mathbf{A}^{\top}\diag(\mathbf{g}) \mathbf{A}$. It is not hard to verify that $\mathbf{B}_{\mathcal{C} \cup \mathcal{M}, \mathcal{C} \cup \mathcal{O}}=\mathbf{A}_{\mathcal{N},\mathcal{C} \cup \mathcal{M}}^{\top}\diag(\mathbf{b}) \mathbf{A}_{\mathcal{N},\mathcal{C} \cup \mathcal{O}}$. Expressing the other submatrices appearing in \eqref{eq:JAt} in a similar fashion and exploiting the invertibility of $\mathbf{A}$ leads to
\begin{equation}\label{eq:Jtilde}
	\tilde{\mathbf{J}}_{A} = \mathbf{A}_{\mathcal{N},\mathcal{C}\cup\mathcal{M}}^{\top}\left(\diag(\mathbf{b})+\diag^2({\mathbf{g}})\diag^{-1}({\mathbf{b}})\right) \mathbf{A}_{\mathcal{N},\mathcal{C}\cup\mathcal{O}}.
\end{equation}

Notice that the middle matrix on the right-hand side of \eqref{eq:Jtilde} is a diagonal matrix $\diag(\tilde{\mathbf{b}})$, where the $n$-th entry of vector $\tilde{\mathbf{b}}$ is $b_n+\frac{g_n^2}{b_n}$ for $n\in \mathcal{N}$. The submatrices of $\mathbf{A}$ appearing in \eqref{eq:Jtilde} can be partitioned column-wise as $\mathbf{A}_{\mathcal{N},\mathcal{C}\cup\mathcal{M}}=[\mathbf{A}_{\mathcal{N},\mathcal{C}}~\mathbf{A}_{\mathcal{N},\mathcal{M}}]$ and $\mathbf{A}_{\mathcal{N},\mathcal{C}\cup\mathcal{O}}=[\mathbf{A}_{\mathcal{N},\mathcal{C}}~\mathbf{A}_{\mathcal{N},\mathcal{O}}]$. Then, matrix $\tilde{\mathbf{J}}_{A}$ is written as:
\begin{equation}\label{eq:JAt2}
	\tilde{\mathbf{J}}_{A} =
	\begin{bmatrix}
	\breve{\mathbf{J}}_A & \breve{\mathbf{J}}_B\\
	\breve{\mathbf{J}}_C  &  \breve{\mathbf{J}}_D
	\end{bmatrix}
\end{equation}
with its submatrices defined as
\begin{align*}\label{eq:JAt3}
	\breve{\mathbf{J}}_A&:=\mathbf{A}_{\mathcal{N},\mathcal{C}}^{\top}\diag(\tilde{\mathbf{b}})\mathbf{A}_{\mathcal{N},\mathcal{C}},~&\breve{\mathbf{J}}_B&:=\mathbf{A}_{\mathcal{N},\mathcal{C}}^{\top}\diag(\tilde{\mathbf{b}})\mathbf{A}_{\mathcal{N},\mathcal{O}},\\
	\breve{\mathbf{J}}_C&:=\mathbf{A}_{\mathcal{N},\mathcal{M}}^{\top}\diag(\tilde{\mathbf{b}})\mathbf{A}_{\mathcal{N},\mathcal{C}},~&\breve{\mathbf{J}}_D&:=\mathbf{A}_{\mathcal{N},\mathcal{M}}^{\top}\diag(\tilde{\mathbf{b}})\mathbf{A}_{\mathcal{N},\mathcal{O}}.
\end{align*}
	
Based on the partitioning of \eqref{eq:JAt2}, the invertibility of $\tilde{\mathbf{J}}_{A}$ can be now characterized using Lemma~\ref{le:struct}: 
	
Firstly, associate $\tilde{\mathbf{J}}_{A}$ with a directed graph $\mathcal{G}=(\mathcal{C}\cup \mathcal{M}\cup \mathcal{O},\mathcal{E})$, where nodes in $\mathcal{C}$ are state vertices, nodes in $\mathcal{M}$ are input vertices, and nodes in $\mathcal{O}$ are output vertices. The edges $\mathcal{E}$ are determined by the sparsity pattern of $\tilde{\mathbf{J}}_{A}$ as in Lemma~\ref{le:struct}. 
	
Secondly, because $\breve{\mathbf{J}}_A$ is a symmetric submatrix of the Laplacian $\mathbf{A}^{\top} \diag(\tilde{\mathbf{b}})\mathbf{A}$, it is invertible; see~\cite[Lemma~1]{MLB15}. 
	
Thirdly, if there exists a vertex-disjoint set of paths connecting $\mathcal{O}$ to $\mathcal{M}$ without going through the substation bus, then Lemma~\ref{le:struct} guarantees that the Schur complement $\breve{\mathbf{J}}_D- \breve{\mathbf{J}}_C \breve{\mathbf{J}}^{-1}_A \breve{\mathbf{J}}_B$ and consequently $\tilde{\mathbf{J}}_A$ are full-rank in general. 
\end{IEEEproof}

\begin{IEEEproof}[Proof of Lemma~\ref{le:Schur}]
Based on the partitioning of \eqref{eq:Jvvp}, let us construct a directed graph $\mathcal{G}=(\mathcal{V},\mathcal{E})$ related to matrix $\mathbf{J}(\mathbf{v},\mathbf{v}')$ according to the procedure detailed before Lemma~\ref{le:struct}. Graph $\mathcal{G}$ has $6(N+1)$ nodes partitioned into three mutually exclusive and collectively exhaustive sets named $\mathcal{X}$, $\mathcal{U}$, and $\mathcal{Y}$. Every node in $\mathcal{X}:=\{x_1,\ldots,x_{2N+2}\}$ is related to an entry of vector $\mathbf{v}$; every node in $\mathcal{U}:=\{u_1,\ldots,u_{2N+2}\}$ is related to an entry of vector $\mathbf{v}'$; and every node in $\mathcal{Y}:=\{y_1,\ldots,y_{2N+2}\}$ is related to one of the specifications \eqref{eq:cpf3:q} and \eqref{eq:cpf2:0}--\eqref{eq:cpf2:p} appearing in the bottom half of $\mathbf{s}(\mathbf{v},\mathbf{v}')$ in \eqref{eq:svv}. The directed edges in $\mathcal{E}$ are drawn based on the sparsity pattern of $\mathbf{J}(\mathbf{v},\mathbf{v}')$, capturing the occurrence of a state variable in the derivative of each one of the specification functions. This graph is not the graph of the underlying physical grid.
	
Due to the structure of $\mathbf{J}(\mathbf{v},\mathbf{v}')$, its diagonal blocks $\mathbf{J}_A(\mathbf{v})$ and $\mathbf{J}_D(\mathbf{v}')$ exhibit the same sparsity pattern for generic values of $(\mathbf{v},\mathbf{v}')$, after row permutations. Since Lemma~\ref{le:JA} guarantees the generic invertibility of $\mathbf{J}_A(\mathbf{v})$ under Criterion~\ref{cr:main}, the same follows for $\mathbf{J}_D(\mathbf{v}')$. By Lemma~\ref{le:bipartite}, the invertibility of $\mathbf{J}_D(\mathbf{v}')$ implies that there exists a perfect matching between the node sets $\mathcal{U}$ and $\mathcal{Y}$. This perfect matching can serve as a set of vertex-disjoint paths from $\mathcal{U}$ and $\mathcal{Y}$; the latter paths are actually single edges and apparently do not pass through any nodes in $\mathcal{X}$. By Lemma~\ref{le:struct}, the invertibility of $\mathbf{J}_A(\mathbf{v})$ and the existence of a set of vertex-disjoint paths from $\mathcal{U}$ to $\mathcal{Y}$ implies that the Schur complement $\mathbf{J}_{D}(\mathbf{v}')-\mathbf{J}_{C}(\mathbf{v})\mathbf{J}_{A}^{-1}(\mathbf{v})\mathbf{J}_{B}(\mathbf{v}')$ is full-rank in general.
\end{IEEEproof}

\bibliographystyle{IEEEtran}
\bibliography{myabrv,power}
\end{document}